# Construction of curve pairs and their applications

**Mehmet Önder**
*Delibekirli Village, Tepe Street, 31440, Kırıkhan, Hatay, Turkey.*
*E-mail: mehmetonder197999@gmail.com*

**Abstract**
In this study, we introduce a new approach to curve pairs by using integral curves. We consider the direction curve and donor curve to study curve couples such as involute- evolute curves, Mannheim partner curves and Bertrand partner curves. We obtain new methods to construct partner curves of a unit speed curve and give some applications related to helices, slant helices and plane curves.

## 1. Introduction

Characterizations of space curve is the most studied subject of the curve theory. These characterizations allows to classify the curves by some properties and are made by two ways; considering a single curve and considering two curves together. Constant slope curves(helices), spherical curves, plane curves and slant helices can be given as the examples of single special curves [1,9,10,16,19]. Particularly, helix curves are used in many applications [2,7,8,15]. Of course, single special curves can be defined by many ways. For instance, if the position vector of a space curve always lies on its rectifying, osculating or normal planes, then the curve is called rectifying curve, osculating curve or normal curve, respectively [4]. These curves also related to some important conditions such as Cesaro's fixed point condition, i.e., a particular point is always involved by Frenet planes of these kind of curves [14]. Especially, there is a connection between Darboux vectors (centrodes) and rectifying curves, which is a useful tool to study on kinematics, mechanics as well as in differential geometry in defining the curves of constant precession [5].

Partner curves are another classifications of curves. These curves are defined by some common properties of two curves. Involute-evolute curves, Bertrand curves, Mannheim curves are the well-known examples of such curves and studied largely [3,11-13,17,18].

Recently, a new approach to curve pairs has been introduced by Choi and Kim [6]. This approach bases on the integral curves of a vector valued function $X$ related to the Frenet frame of a Frenet curve $\alpha$. They have defined principal-donor and principal-directional curves in $E^3$ by this approach and given properties.

In this study, we consider integral curves to study partner curves such as involute and evolute curves, Mannheim partner curves and Bertrand partner curves. We define these curves as direction curves and obtain relationships between their Frenet elements. Furthermore, we apply the obtainer results to some special curves such as helices, plane curves and slant helices.

## 2. Preliminaries

A unit speed curve $\alpha : I \to E^3$ is called a helix or general helix if its unit tangent vector $T(s) = \alpha'(s)$ always makes a constant angle $\varphi \neq \pi/2$ with a constant vector $\ell$, i.e., $\langle T, \ell \rangle = \cos\varphi$ is constant for all $s \in I$. The function $\kappa(s) = \|\alpha''(s)\|$ is called curvature (or first curvature) of $\alpha$ and $\alpha$ is called Frenet curve, if $\kappa(s) \neq 0$. The vector given $N$ by the equality $\alpha''(s) = \kappa(s) N(s)$ is called unit principal normal vector of $\alpha$. Then, we define the unit vector $B(s)=T(s)\times N(s)$ as the unit binormal vector of $\alpha$. Thus, $\{T,N,B\}$ is an orthonormal frame along $\alpha$ and known as the Frenet frame of $\alpha$ which has the derivative formulas

$$\begin{bmatrix} T' \\ N' \\ B' \end{bmatrix} = \begin{bmatrix} 0 & \kappa & 0 \\ -\kappa & 0 & \tau \\ 0 & -\tau & 0 \end{bmatrix} \begin{bmatrix} T \\ N \\ B \end{bmatrix} \tag{1}$$

where $\tau(s)$ is the second curvature (or torsion) of $\alpha$ at $s$. Special curves are characterized by their curvatures. For instance, a curve $\alpha$ is a helix if and only if $\frac{\tau}{\kappa}(s)=\text{constant}$ [16]. A curve $\alpha$ in $E^3$ with non-zero curvature $\kappa(s)$ is called a slant helix if its principal normal vectors make a constant angle with a fixed direction and characterized by the differential equation of its curvature $\kappa$ and its torsion $\tau$ given by

$$\frac{\kappa^2}{\left(\kappa^2+\tau^2\right)^{3/2}}\left(\frac{\tau}{\kappa}\right)' = \text{constant}.$$

(See [10]).

Now, we give the definitions of some associated curves defined by Choi and Kim [6]. Let $I\subset\mathbb{R}$ be an open interval. For a unit speed Frenet curve $\alpha(s): I\rightarrow E^3$, consider a vector valued function $X$ defined by

$$X(s)=a(s)T(s)+b(s)N(s)+c(s)B(s), \tag{2}$$

where $a$, $b$ and $c$ are smooth functions of $s$. Let

$$a^2(s)+b^2(s)+c^2(s)=1, \tag{3}$$

holds. Then the definitions of $X$-donor curve and $X$-direction curve in $E^3$ are given as follows.

**Definition 2.1. ([6])** Let $\alpha$ be a Frenet curve in Euclidean 3-space $E^3$ and $X$ be a unit vector valued function with properties (2) and (3). The curve $\gamma: I\rightarrow E^3$ which is also an integral curve of $X$ is called an $X$-direction curve of $\alpha$. The curve $\alpha$ having $\gamma$ as an $X$-direction curve is called the $X$-donor curve of $\gamma$ in $E^3$.

**Definition 2.2. ([6])** In $E^3$, an integral curve of principal normal vector $N(s)$ of $\alpha$ in (2) is called the principal-direction curve of $\alpha$.

**Remark 2.1. ([6])** A principal-direction curve is an integral curve of $X(s)$ with $a(s)=c(s)=0,\ b(s)=1$ for all $s$ in (2).

Let $\alpha=\alpha(s)$ be a Frenet curve in $E^3$ with arc length parameter $s$ and $X$ be a unit vector field as given in (2) and (3). Let $\beta: I\rightarrow E^3$ be an $X$-direction curve of $\alpha$. The Frenet vectors and curvatures of $\alpha$ and $\beta$ be denoted by $\{T,N,B\}$, $\kappa, \tau$ and $\{\bar{T},\bar{N},\bar{B}\}$, $\bar{\kappa}, \bar{\tau}$, respectively. From (2) and (3) it is clear that arclength parameter $s$ of $\alpha$ can be also taken as arclength parameter of $\beta$. Then, differentiating (2) with respect to $s$ and using the fact that $X=\bar{T}$, it follows

$$\bar{\kappa}\bar{N}=(a'-b\kappa)T+(b'+a\kappa-c\tau)N+(c'+b\tau)B \tag{4}$$

Equality (4) allows us to study curve pairs such as involute and evolute curves, Mannheim partner curves and Bertrand partner curve. So, in the following section, we will assumed that the unit vector $X$ is defined as given in (2) and satisfied (3).

## 3. Involute-evolute-direction curves and their applications

In this section, we define involute-evolute-direction curves and obtain relationships between this curves.

**Definition 3.1.** Let $\alpha : I \to E^3$ be a Frenet curve in Euclidean 3-space $E^3$ and $\beta : I \to E^3$ be $X$-direction curve of $\alpha$. If $\beta$ is an evolute of $\alpha$, then $\beta$ is called evolute-direction curve of $\alpha$. Then, $\alpha$ is said to be involute-donor curve of $\beta$.

**Proposition 3.2.** *For the Frenet curve $\alpha : I \to E^3$, the curve $\beta$ is evolute-direction curve of $\alpha$ if and only if*

$$a(s) = 0\,,\ \ b(s) = \sin\left(\int \tau ds\right),\ c(s) = \cos\left(\int \tau ds\right). \tag{5}$$

**Proof.** From the definition of involute-evolute curves, it is well-known that $\langle \bar{T}, T \rangle = 0$ and $\bar{N} = T$. Then, from (4) we have that $\beta$ is evolute-direction curve of $\alpha$ if and only if $a = 0$ and the system

$$\begin{cases} -b\kappa = \bar{\kappa} \neq 0, \\ b' - c\tau = 0, \\ c' + b\tau = 0. \end{cases} \tag{6}$$

holds, The solution of system (6) is $\left\{ b = \sin\left(\int \tau ds\right),\ c = \cos\left(\int \tau ds\right) \right\}$.

From Proposition 3.2, we can give the following definition.

**Definition 3.3.** *An integral curve of vector field* $\sin\left(\int \tau ds\right) N + \cos\left(\int \tau ds\right) B$ *is called evolute-direction curve of $\alpha$.*

Proposition 3.2 gives a method to construct an evolute of a curve by using its Frenet elements $N, B$ and torsion $\tau$. It means that to obtain an evolute curve of a curve, it is enough to know the Frenet elements $N, B, \tau$ of reference curve.

In the following characterizations, we give the relationships between curvatures and Frenet vectors of involute-evolute-direction curves.

**Corollary 3.4.** *The relations between the Frenet vectors of involute-evolute-direction curves are given by*

$$\bar{T}(s) = \sin\left(\int \tau ds\right) N(s) + \cos\left(\int \tau ds\right) B(s)\,,\ \ \bar{N}(s) = T(s) \tag{7}$$

$$\bar{B}(s) = \cos\left(\int \tau ds\right) N(s) - \sin\left(\int \tau ds\right) B(s)\,. \tag{8}$$

**Proof.** The proof is clear from (5) and definition of involute-evolute curves.

**Theorem 3.5.** *Let $\alpha$ be a Frenet curve in $E^3$ and $\beta$ be evolute-direction curve of $\alpha$. Then*

$$\bar{\kappa} = \kappa \left| \sin\left(\int \tau ds\right) \right|,\ \ \bar{\tau} = \kappa \cos\left(\int \tau ds\right), \tag{9}$$

$$\kappa=\sqrt{\bar{\kappa}^2+\bar{\tau}^2}\,,\ \ \tau=\frac{\bar{\kappa}^2}{\bar{\kappa}^2+\bar{\tau}^2}\left(\frac{\bar{\tau}}{\bar{\kappa}}\right)'. \tag{10}$$

**Proof.** From (5) and the first equation of system (6), we have the first equality of (9) immediately. Moreover from (8), it follows $\bar{B}'=-\kappa\cos\left(\int\tau ds\right)T$. Since $\bar{N}=T$, we have

$$\bar{\tau}=-\left\langle \bar{B}',\bar{N}\right\rangle=\kappa\cos\left(\int\tau ds\right). \tag{11}$$

Let now obtain the relations given in (10). Differentiating $\bar{N}=T$ and taking the norm of the obtained equality gives $\kappa=\sqrt{\bar{\kappa}^2+\bar{\tau}^2}$. For the second relation, considering (11) and $\kappa=\sqrt{\bar{\kappa}^2+\bar{\tau}^2}$, we have

$$\int\tau ds=\arccos\left(\frac{\bar{\tau}}{\sqrt{\bar{\kappa}^2+\bar{\tau}^2}}\right). \tag{12}$$

Differentiating (12), it follows $\tau=\dfrac{\bar{\kappa}^2}{\bar{\kappa}^2+\bar{\tau}^2}\left(\dfrac{\bar{\tau}}{\bar{\kappa}}\right)'$.

From (10) we have the following important corollary which allow us to apply the notion of involute-evolute-direction curves to helices and slant helices.

**Corollary 3.6.** Let $\beta$ *be evolute-direction curve of* $\alpha$. *Then*

$$\frac{\tau}{\kappa}=\frac{\bar{\kappa}^2}{\left(\bar{\kappa}^2+\bar{\tau}^2\right)^{3/2}}\left(\frac{\bar{\tau}}{\bar{\kappa}}\right)' \tag{13}$$

*holds.*

### 3.1. Applications of involute-evolute-direction curves

In this section, we focus on relations between involute-evolute-direction curves and other important curves such as general helices, slant helices and plane curves in $E^3$.

Considering Corollary 3.6, we have the following theorems which gives a way to construct the examples of slant helices by using general helices.

***Theorem 3.7.*** *For the Frenet curve* $\alpha:I\to E^3$ *with an evolute-direction curve* $\beta$, *the followings statements are equivalent,*

*i)* $\alpha$ *is a helix curve.*

*ii)* $\alpha$ *is an involute-donor curve of a slant helix.*

*iii) An evolute-direction curve of* $\alpha$ *is a slant helix.*

***Theorem 3.8.*** *For the Frenet curve* $\alpha:I\to E^3$ *with an evolute-direction curve* $\beta$, *the followings statements are equivalent,*

*i)* $\alpha$ *is a plane curve.*

*ii)* $\alpha$ *is an involute-donor curve of a helix.*

*iii) An evolute-direction curve of* $\alpha$ *is a helix.*

**Example 3.9.** Let consider the general helix given by the parametrization $\alpha(s)=\left(\cos\frac{s}{\sqrt{2}},\ \sin\frac{s}{\sqrt{2}},\ \frac{s}{\sqrt{2}}\right)$ in $E^3$ (Fig. 1). The Frenet elements of $\alpha$ are

$$T(s)=\left(-\frac{1}{\sqrt{2}}\sin\frac{s}{\sqrt{2}},\ \frac{1}{\sqrt{2}}\cos\frac{s}{\sqrt{2}},\ \frac{1}{\sqrt{2}}\right),$$

$$N(s)=\left(-\cos\frac{s}{\sqrt{2}},\ \sin\frac{s}{\sqrt{2}},\ 0\right),$$

$$B(s)=\left(\frac{1}{\sqrt{2}}\sin\frac{s}{\sqrt{2}},\ -\frac{1}{\sqrt{2}}\cos\frac{s}{\sqrt{2}},\ \frac{1}{\sqrt{2}}\right),$$

$$\kappa=\tau=\frac{1}{2}.$$

From (6) we have, $X(s)=\left(x_1(s),\ x_2(s),\ x_3(s)\right)$ where

$$x_1(s)=-\sin\left(\frac{s}{2}+c\right)\cos\frac{s}{\sqrt{2}}+\frac{1}{\sqrt{2}}\cos\left(\frac{s}{2}+c\right)\sin\frac{s}{\sqrt{2}},$$

$$x_2(s)=\sin\left(\frac{s}{2}+c\right)\sin\frac{s}{\sqrt{2}}-\frac{1}{\sqrt{2}}\cos\left(\frac{s}{2}+c\right)\cos\frac{s}{\sqrt{2}},$$

$$x_3(s)=\frac{1}{\sqrt{2}}\cos\left(\frac{s}{2}+c\right).$$

and $c$ is integration constant. Now, we can construct a slant helix $\beta$ which is also an evolute-direction curve of $\alpha$:

$$\beta=\int_0^s\beta'(s)ds=\int_0^s X(s)ds=\left(\beta_1(s),\ \beta_2(s),\beta_3(s)\right),$$

where

$$\beta_1(s)=\int_0^s\left[-\sin\left(\frac{s}{2}+c\right)\cos\frac{s}{\sqrt{2}}+\frac{1}{\sqrt{2}}\cos\left(\frac{s}{2}+c\right)\sin\frac{s}{\sqrt{2}}\right]ds,$$

$$\beta_2(s)=\int_0^s\left[\sin\left(\frac{s}{2}+c\right)\sin\frac{s}{\sqrt{2}}-\frac{1}{\sqrt{2}}\cos\left(\frac{s}{2}+c\right)\cos\frac{s}{\sqrt{2}}\right]ds,$$

$$\beta_3(s)=\int_0^s\frac{1}{\sqrt{2}}\cos\left(\frac{s}{2}+c\right)ds.$$

By taking $c=0$, we have Fig. 2

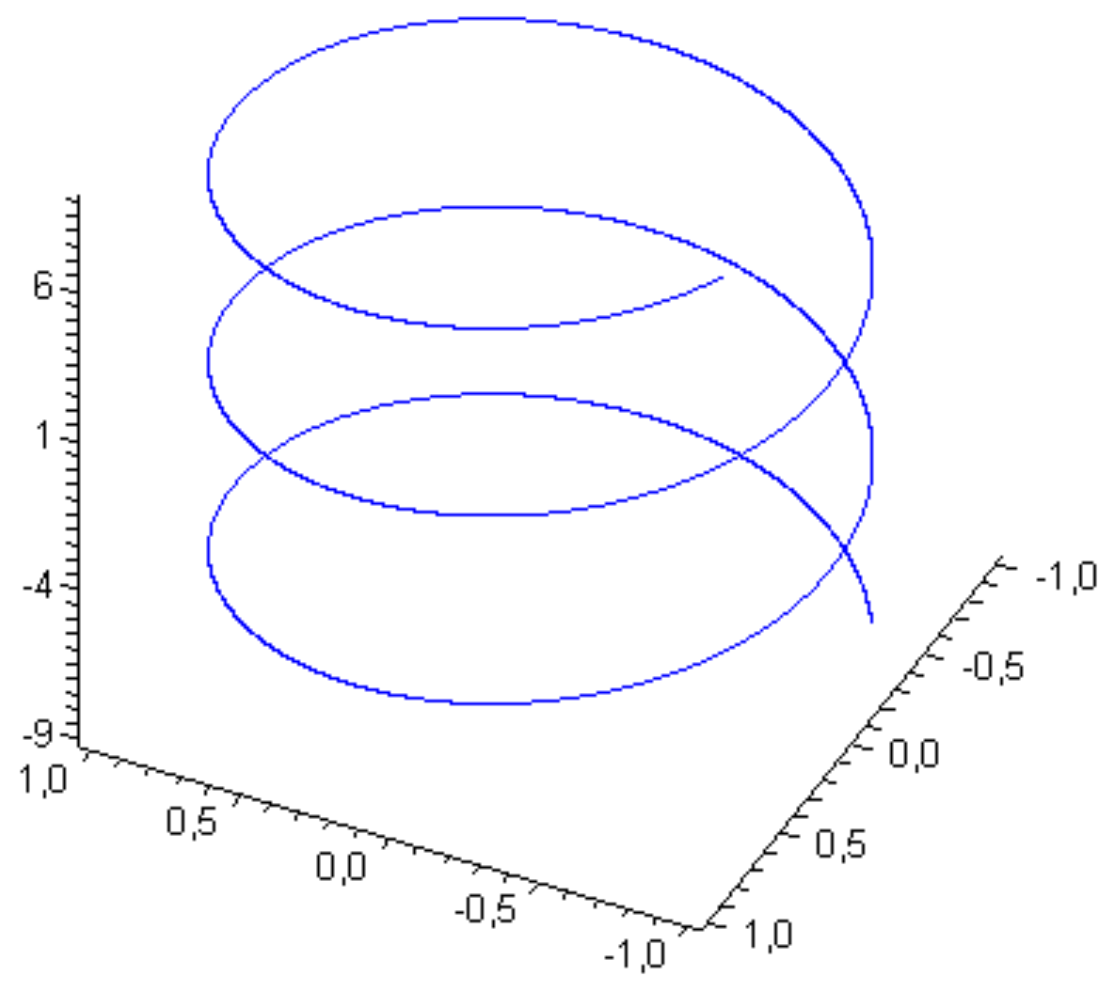

**Fig. 1.** General helix $\alpha$.

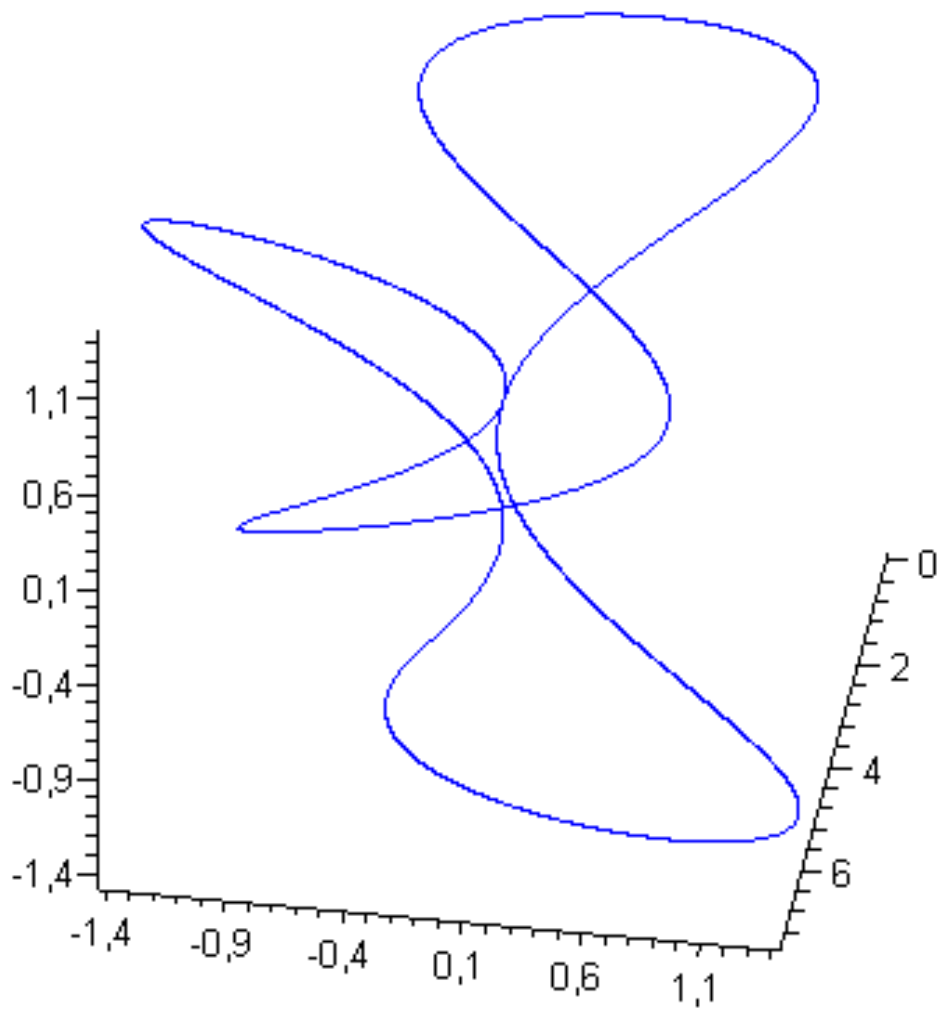

**Fig. 2.** Slant helix $\beta$ constructed by $\alpha$.

**4. Mannheim-direction curves and their applications**

In this section, we define Mannheim-direction curves and obtain relationships between this curves.

**Definition 4.1.** Let $\alpha : I \to E^3$ be a Frenet curve in $E^3$ and $\beta : I \to E^3$ be $X$-direction curve of $\alpha$. If $\beta$ is a Mannheim curve of $\alpha$ and $\alpha$ is a Mannheim partner curve of $\beta$, then $\beta$ is called Mannheim-direction curve of $\alpha$. Then, $\alpha$ is said to be Mannheim-donor curve of $\beta$.

**Proposition 4.2.** *Let $\alpha$ be a Frenet curve in $E^3$. Then $\beta$ is Mannheim-direction curve of $\alpha$ if and only if*

$$a(s) = \sin\left(\int \kappa ds\right),\ \ b(s) = \cos\left(\int \kappa ds\right),\ \ c(s) = 0. \tag{14}$$

**Proof.** From the definition of Mannheim curves, it is well-known that $\bar{N} = B$. Then, from (4) we have that $\beta$ is a Mannheim-direction curve of $\alpha$ if and only if the system,

$$\begin{cases} a' - b\kappa = 0, \\ b' + a\kappa - c\tau = 0, \\ c' + b\tau = \bar{\kappa} \neq 0. \end{cases} \tag{15}$$

holds. Multiplying the first equation in (15) with $a$ and second equation with $b$ and adding the results gives

$$aa' + bb' - bc\tau = 0. \tag{16}$$

From (3) and (16) we have $c(c' + b\tau) = 0$. Since $c' + b\tau = \bar{\kappa} \neq 0$ in (15), it follows $c = 0$. Then the system (15) reduced to the following system,

$$\begin{cases} a' - b\kappa = 0, \\ b' + a\kappa = 0, \\ b\tau = \bar{\kappa} \neq 0. \end{cases} \tag{16}$$

Then the solution of system (15) is obtained as $\left\{a = \sin\left(\int \kappa ds\right),\ b = \cos\left(\int \kappa ds\right),\ c = 0\right\}$.

**Definition 4.3.** *An integral curve of vector field $\sin\left(\int \kappa ds\right)T + \cos\left(\int \kappa ds\right)N$ is called Mannheim-direction curve of $\alpha$.*

Proposition 4.2 gives a method to construct a Mannheim curve of a unit speed curve by using its Frenet elements $T, N$ and $\kappa$. It means that to obtain a Mannheim curve of a curve, it is enough to know the Frenet elements $T, N, \kappa$ of reference curve. Moreover, we have that a Mannheim-direction curve of a space curve is also osculating direction curve of the same curve.

In the following characterizations, we give the relationships between curvatures and Frenet vectors of Mannheim-direction curves and their applications.

**Corollary 4.4.** *The relations between the Frenet vectors of Mannheim-direction curves are given by*

$$\bar{T}(s) = \sin\left(\int \kappa ds\right) T(s) + \cos\left(\int \kappa ds\right) N(s), \quad \bar{N}(s) = B(s) \tag{17}$$

$$\bar{B}(s) = \cos\left(\int \kappa ds\right) T(s) - \sin\left(\int \kappa ds\right) N(s). \tag{18}$$

**Proof.** The proof is clear from the obtained results.

**Theorem 4.5.** *Let $\alpha$ be a Frenet curve in $E^3$ and $\beta$ be Mannheim-direction curve of $\alpha$. Then*

$$\bar{\kappa} = \left| \tau \cos\left(\int \kappa ds\right) \right|, \quad \bar{\tau} = \tau \sin\left(\int \kappa ds\right), \tag{19}$$

$$\kappa = \left| \frac{\bar{\kappa}^2}{\bar{\kappa}^2 + \bar{\tau}^2} \left(\frac{\bar{\tau}}{\bar{\kappa}}\right)' \right|, \quad |\tau| = \sqrt{\bar{\kappa}^2 + \bar{\tau}^2}. \tag{20}$$

**Proof.** From (16) and the second equation of (14), we have the first equality of (19) immediately. Moreover from (18), it follows $\bar{B}' = -\tau \sin\left(\int \kappa ds\right) B$. Since $\bar{N} = B$, we have

$$\bar{\tau} = -\langle \bar{B}', \bar{N} \rangle = \tau \sin\left(\int \kappa ds\right). \tag{21}$$

Let now obtain the relations given in (20). Differentiating the equality $\bar{N} = B$ and taking the norm of the obtained result gives $|\tau| = \sqrt{\bar{\kappa}^2 + \bar{\tau}^2}$. For the first relation in (20), considering (21) and $|\tau| = \sqrt{\bar{\kappa}^2 + \bar{\tau}^2}$, we have

$$\int \kappa ds = \arccos\left( \pm \frac{\bar{\kappa}}{\sqrt{\bar{\kappa}^2 + \bar{\tau}^2}} \right). \tag{22}$$

Differentiating (22), it follows $\kappa = \left| \frac{\bar{\kappa}^2}{\bar{\kappa}^2 + \bar{\tau}^2} \left(\frac{\bar{\tau}}{\bar{\kappa}}\right)' \right|$.

From (20) we have the following important corollary.

**Corollary 4.6.** Let *$\beta$ be a Mannheim-direction curve of $\alpha$. Then*

$$\frac{\tau}{\kappa} = \pm \frac{1}{\left| \frac{\bar{\kappa}^2}{\left(\bar{\kappa}^2 + \bar{\tau}^2\right)^{3/2}} \left(\frac{\bar{\tau}}{\bar{\kappa}}\right)' \right|} \tag{23}$$

*holds.*

### 4.1. Applications of Mannheim-direction curves

In this section, we focus on relations between Mannheim-direction curves, helices and slant helices in $E^3$.

Considering Corollary 4.6, we have the following theorems which gives a way to construct the examples of slant helices by using general helices for which the obtained slant helix is a Mannheim curve of the helix.

***Theorem 4.7.*** *Let $\alpha : I \to E^3$ be a Frenet curve in $E^3$ and $\beta$ be a Mannheim-direction curve of $\alpha$. Then, we have that the statements*

*i) $\alpha$ is a helix.*

*ii) $\alpha$ is a Mannheim-donor curve of a slant helix.*

*iii) A Mannheim-direction curve of $\alpha$ is a slant helix.*

*are equivalent.*

**Example 4.8.** Let consider the general helix given in Example 3.9. From (14) we have, $X(s) = \left(x_1(s),\, x_2(s),\, x_3(s)\right)$ where

$$x_1(s) = -\frac{1}{\sqrt{2}}\sin\left(\frac{s}{2}+c\right)\sin\left(\frac{s}{\sqrt{2}}\right) - \cos\left(\frac{s}{2}+c\right)\cos\left(\frac{s}{\sqrt{2}}\right),$$

$$x_2(s) = \frac{1}{\sqrt{2}}\sin\left(\frac{s}{2}+c\right)\cos\left(\frac{s}{\sqrt{2}}\right) + \cos\left(\frac{s}{2}+c\right)\sin\left(\frac{s}{\sqrt{2}}\right),$$

$$x_3(s) = \frac{1}{\sqrt{2}}\sin\left(\frac{s}{2}+c\right),$$

and $c$ is integration constant. Now, we can construct a slant helix $\beta$ which is also a Mannheim-direction of $\alpha$ by taking

$$\beta = \int_0^s \beta'(s)ds = \int_0^s X(s)ds = \left(\beta_1(s),\, \beta_2(s), \beta_3(s)\right).$$

Choosing $c = 0$, we have the following figures representing $\beta$.

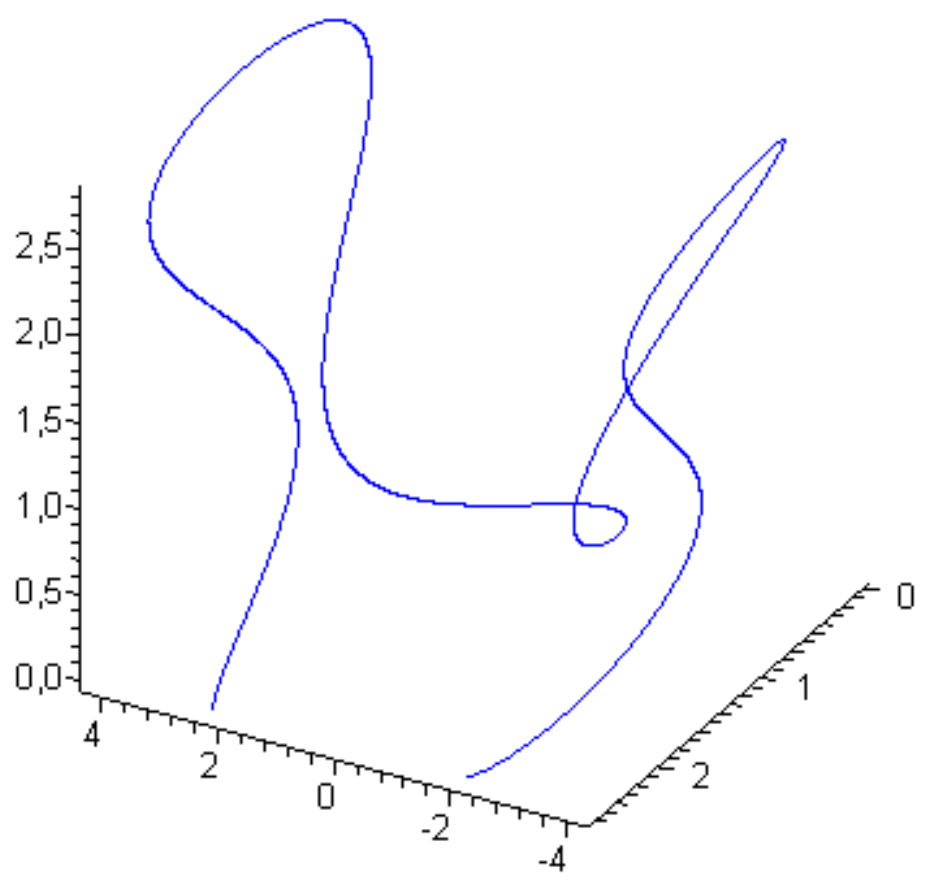

**Fig 3.** The curve $\beta$ for $-4\pi \le s \le 4\pi$.

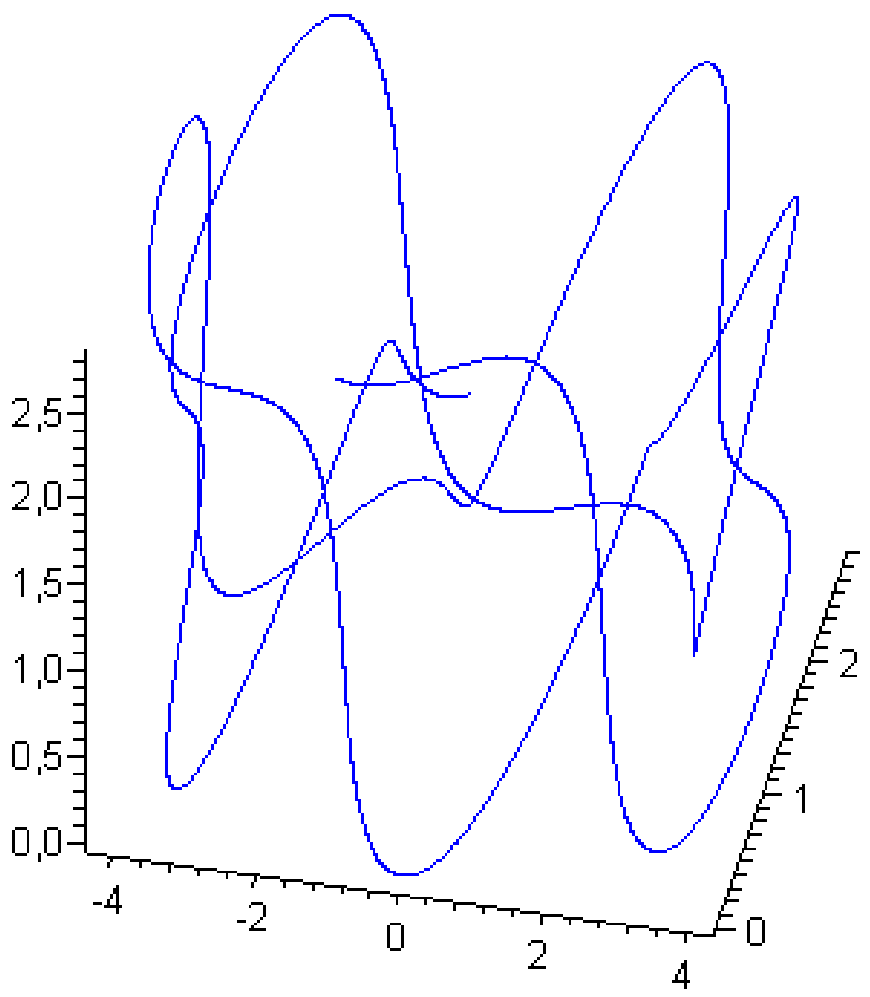

**Fig 4.** The curve $\beta$ for $-10\pi \le s \le 10\pi$.

## 5. Bertrand-direction curves

In this section, we define Bertrand-direction curves and obtain relationships between this curves.

**Definition 5.1.** Let $\alpha : I \to E^3$ be a Frenet curve in Euclidean 3-space $E^3$ and $\beta : I \to E^3$ be $X$ -direction curve of $\alpha$. If $\beta$ is a Bertrand curve of $\alpha$ and $\alpha$ is a Bertrand partner curve of $\beta$, then $\beta$ is called Bertrand-direction curve of $\alpha$. Then, $\alpha$ is said to be Bertrand-donor curve of $\beta$.

**Proposition 5.2.** *Let $\alpha$ be a Frenet curve in Euclidean 3-space $E^3$. Then $\beta$ is Bertrand-direction curve of $\alpha$ if and only if*

$$a = \cos\theta,\ b = 0,\ c = \sin\theta, \tag{24}$$

*where $\theta$ is the constant angle between the tangent lines of the curves.*

**Proof.** From the definition of Bertrand curves, it is well-known that $\bar{N} = N$. Then, from (4) we have that $\beta$ is a Bertrand-direction curve of $\alpha$ if and only if the system

$$\begin{cases} a' - b\kappa = 0, \\ b' + a\kappa - c\tau = \bar{\kappa} \neq 0, \\ c' + b\tau = 0. \end{cases} \tag{25}$$

holds. Multiplying the first, the second and the third equations in (25) with $a$, $b$ and $c$, respectively, and adding the results gives

$$aa' + bb' + cc' = b\bar{\kappa}. \tag{26}$$

From (3) and (26) we have $b = 0$. Then the solution of system (25) is obtained as $\{a = \text{constant},\ b = 0,\ c = \text{constant}\}$ and we have $X = aT + cB$. Since $\beta$ is the integral curve of $X$, if the constant angle between the tangent lines of $\alpha$ and $\beta$ is taken as $\theta$, then the solution of the system (25) can be given as $\{a = \cos\theta,\ b = 0,\ c = \sin\theta;\ \theta\text{=constant}\}$

**Definition 5.3.** *An integral curve of vector field $\cos\theta T + \sin\theta B$, where $\theta = \text{constant}$, is called Bertrand-direction curve of $\alpha$.*

Proposition 5.2 gives a method to construct a Bertrand curve of a unit speed curve by using its Frenet vectors $T, B$ and a constant $\theta$. It means that to obtain a Bertrand curve of a curve, it is enough to know the Frenet elements $T, B$ of reference curve and a constant. Moreover, we have that a Bertrand-direction curve of a space curve is also rectifying-direction curve of the same curve.

In the following characterizations, we give the relationships between curvatures and Frenet vectors of Bertrand-direction curves and their applications.

**Corollary 5.4.** *The relations between the Frenet vectors of Bertrand-direction curves are given by*

$$\bar{T}(s) = \cos\theta T(s) + \sin\theta B(s),\ \bar{N}(s) = N(s),\ \bar{B}(s) = -\sin\theta T(s) + \cos\theta B(s). \tag{27}$$

**Proof.** The proof is clear from the obtained results.

**Theorem 5.5.** *Let $\alpha$ be a Frenet curve in $E^3$ and $\beta$ be Bertrand-direction curve of $\alpha$. Then*

$$\bar{\kappa} = \left|\kappa\cos\theta - \tau\sin\theta\right|,\ \bar{\tau} = \kappa\sin\theta + \tau\cos\theta, \tag{28}$$

$$\kappa = \left|\bar{\kappa}\cos\theta + \bar{\tau}\sin\theta\right|,\ \ \tau = \bar{\tau}\cos\theta - \bar{\kappa}\sin\theta . \tag{29}$$

**Proof.** Differentiating the first equality in (27) gives $\bar{\kappa}\bar{N} = (\kappa\cos\theta - \tau\sin\theta)N$ which gives the first equality in (28) immediately. Moreover from the third equality of (27), it follows $\bar{B}' = -(\kappa\sin\theta + \tau\cos\theta)N$. Since $\bar{N} = N$, we have

$$\bar{\tau} = -\left\langle \bar{B}', \bar{N} \right\rangle = \kappa\sin\theta + \tau\cos\theta . \tag{30}$$

Let now obtain the relations given in (29). From (27), we can write

$$T(s) = \cos\theta\bar{T}(s) - \sin\theta\bar{B}(s),\ \ B(s) = \sin\theta\bar{T}(s) + \cos\theta\bar{B}(s). \tag{31}$$

Differentiating this equalities and using the similar way we have (29).

From Theorem 5.5, we have the following important corollary.

**Corollary 5.6.** Let $\beta$ *be a Bertrand-direction curve of* $\alpha$. *Then*

$$\frac{\tau}{\kappa} = \frac{-\sin\theta + \left(\bar{\tau}/\bar{\kappa}\right)\cos\theta}{\left|\cos\theta + \left(\bar{\tau}/\bar{\kappa}\right)\sin\theta\right|} \tag{32}$$

$$\frac{\kappa^2}{\left(\kappa^2+\tau^2\right)^{3/2}}\left(\frac{\tau}{\kappa}\right)' = \frac{\bar{\kappa}^2}{\left(\bar{\kappa}^2+\bar{\tau}^2\right)^{3/2}}\left(\frac{\bar{\tau}}{\bar{\kappa}}\right)' \tag{33}$$

*hold.*

### 5.1. Applications of Bertrand-direction curves

In this section, by considering the obtained results, we give relations between Bertrand-direction curves and some other curves such as helices, plane curves and slant helices in $E^3$.

Considering Corollary 5.6, we have the following theorems which gives a way to construct the examples of helices and slant helices for which the obtained curve is a Bertrand curve of the reference curve.

***Theorem 5.7.*** *Let* $\alpha : I \to E^3$ *be a Frenet curve in* $E^3$ *and* $\beta$ *be a Bertrand-direction curve of* $\alpha$. *Then the statements*

*i)* $\alpha$ *is a helix in* $E^3$.

*ii)* $\alpha$ *is a Bertrand-donor curve of a helix.*

*iii) A Bertrand-direction curve of* $\alpha$ *is a helix.*

*are equivalent.*

***Theorem 5.8.*** *Let* $\alpha : I \to E^3$ *be a Frenet curve in* $E^3$ *and* $\beta$ *be a Bertrand-direction curve of* $\alpha$. *If* $\alpha$ *is a plane curve, then* $\beta$ *is a helix. Similarly, if* $\beta$ *is a plane curve, then* $\alpha$ *is a helix.*

***Theorem 5.9.*** *Let* $\alpha : I \to E^3$ *be a Frenet curve in* $E^3$ *and* $\beta$ *be a Bertrand-direction curve of* $\alpha$. *Then the followings are equivalent,*

*i)* $\alpha$ *is a slant helix.*

*ii)* $\alpha$ *is a Bertrand-donor curve of a slant helix.*

*iii) A Bertrand-direction curve of* $\alpha$ *is a slant helix.*

These theorems allow to construct Bertrand curves which are also special curves. Theorem 5.7 says that a Bertrand curve of a helix is also a helix. Similarly, Theorem 5.9 says that a Bertrand curve of slant helix is also a slant helix. The result of Theorem 5.8 is classical that a helix can be constructed from a plane curve []. As applications of theorems, we can consider the following examples.

**Example 5.10.** Let consider the general helix given in Example 3.9. By choosing $\theta = \pi/3$, from (24) we have,

$$X(s) = \left( \frac{\sqrt{3}-1}{2\sqrt{2}} \sin\left(\frac{s}{\sqrt{2}}\right), \frac{1-\sqrt{3}}{2\sqrt{2}} \cos\left(\frac{s}{\sqrt{2}}\right), \frac{1+\sqrt{3}}{2\sqrt{2}} \right) \tag{34}$$

Now, we can construct a helix $\beta$ which is also a Bertrand-direction curve of $\alpha$ by taking

$$\begin{aligned} \beta &= \int_0^s \beta'(s)ds \\ &= \int_0^s X(s)ds = \left( -\frac{\sqrt{3}-1}{2} \cos\left(\frac{s}{\sqrt{2}}\right), \frac{1-\sqrt{3}}{2} \sin\left(\frac{s}{\sqrt{2}}\right), \frac{1+\sqrt{3}}{2\sqrt{2}} s \right) \end{aligned}$$

(Fig 5)

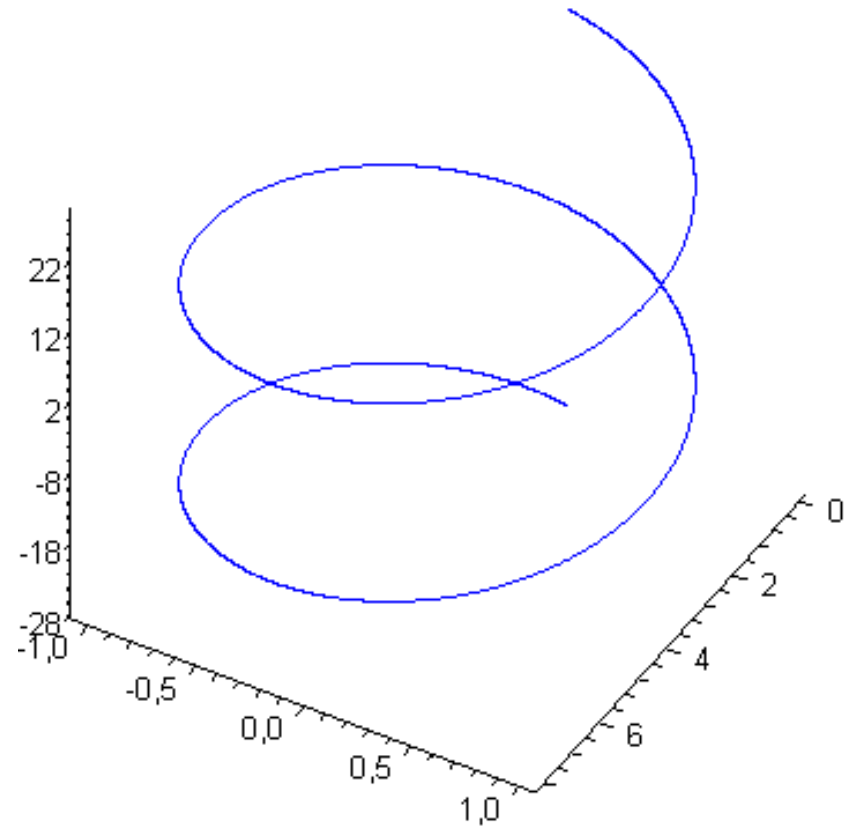

**Fig 5.** The Bertrand direction curve $\beta$ for $-4\pi \le s \le 4\pi$.

**Example 5.11.** Let consider the slant helix $\alpha$ given by the parametrization,

$$\begin{aligned} \alpha(s) = \frac{9}{\sqrt{82}} \Bigg( & \frac{\sqrt{82}-82}{8(41+\sqrt{82})} \sin\left(\left(1+\frac{\sqrt{82}}{41}\right)s\right) + \frac{\sqrt{82}+82}{8(\sqrt{82}-41)} \sin\left(\left(1-\frac{\sqrt{82}}{41}\right)s\right) - \frac{1}{2}\sin s, \\ & \frac{82-\sqrt{82}}{8(41+\sqrt{82})} \cos\left(\left(1+\frac{\sqrt{82}}{41}\right)s\right) - \frac{\sqrt{82}+82}{8(\sqrt{82}-41)} \cos\left(\left(1-\frac{\sqrt{82}}{41}\right)s\right) + \frac{1}{2}\cos s, \\ & \frac{9}{4}\cos\left(\frac{\sqrt{82}}{41}s\right) \Bigg). \end{aligned}$$

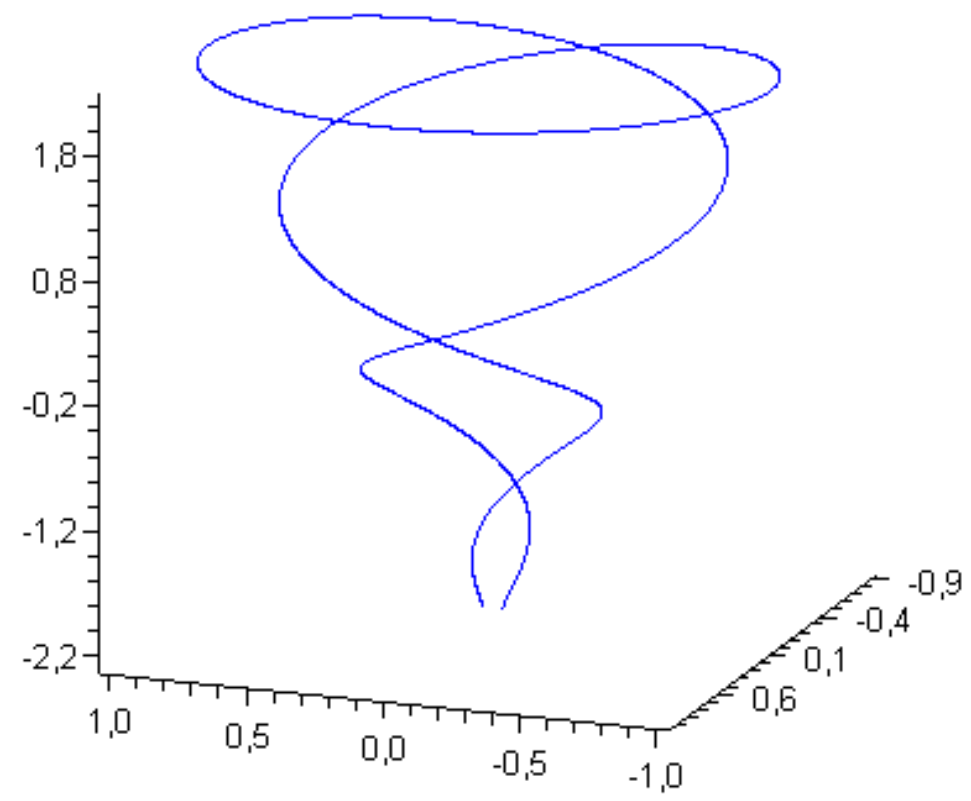

**Fig. 6.** Slant helix $\alpha$.

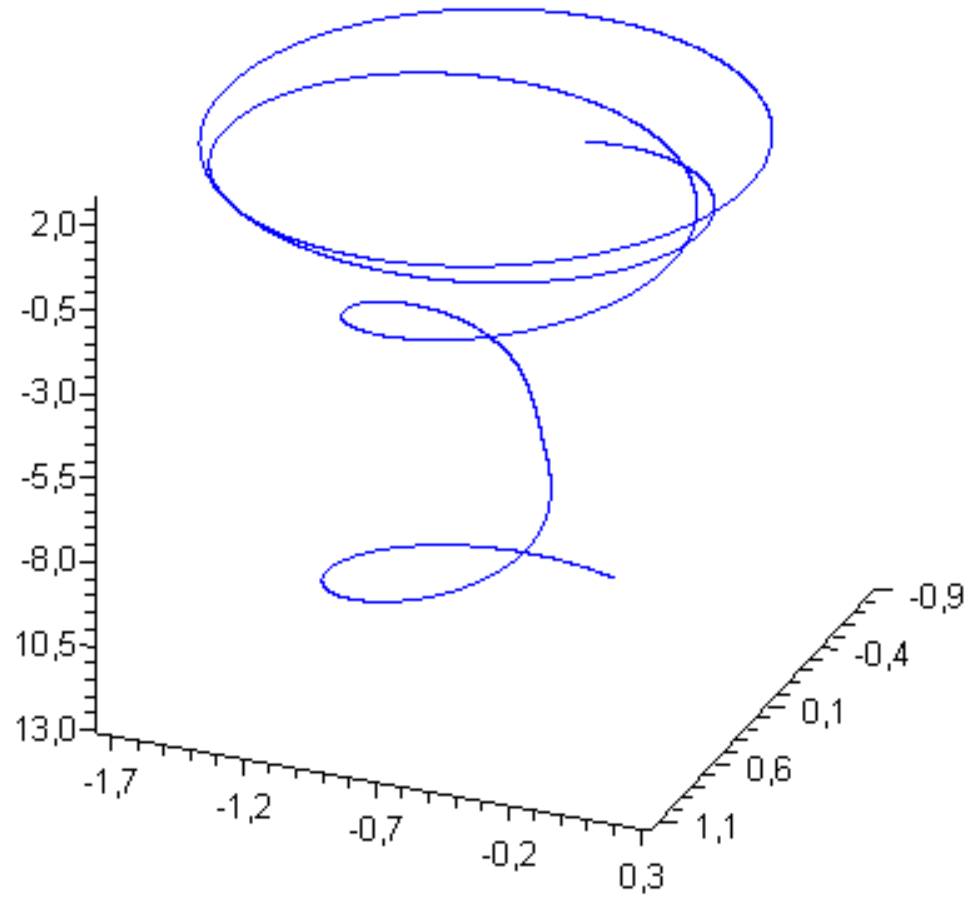

**Fig. 7.** Bertrand direction curve $\beta$.

By computing Frenet frames, choosing $\theta = \pi/4$ and using the fact that we obtain another slant helix given by the parametrization $\beta = \int_0^s X(s)ds = (\beta_1, \beta_2, \beta_3)$ which is also a Bertrand-direction curve of $\alpha$, where

$$\beta_1(s) = \frac{\sqrt{2}}{2}\int_0^s \left[ \cos\left(\frac{\sqrt{82}}{82}s\right)\left(\frac{\sqrt{82}}{82}\sin s - \cos s\right) - \sin\left(\frac{\sqrt{82}}{82}s\right)\left(\frac{\sqrt{82}}{82}\sin s + \cos s\right)\right] ds,$$

$$\beta_2(s) = \frac{\sqrt{2}}{2}\int_0^s \left[ -\cos\left(\frac{\sqrt{82}}{82}s\right)\left(\frac{\sqrt{82}}{82}\cos s + \sin s\right) + \sin\left(\frac{\sqrt{82}}{82}s\right)\left(\frac{\sqrt{82}}{82}\cos s - \sin s\right)\right] ds$$

$$\beta_3(s) = \frac{9\sqrt{164}}{164}\int_0^s \left[ \cos\left(\frac{\sqrt{82}}{82}s\right) - \sin\left(\frac{\sqrt{82}}{82}s\right)\right] ds$$

(Fig. 7).